\documentclass{article} 
\usepackage{times}
\usepackage{graphicx}
\usepackage{amssymb}
\usepackage{url,hyperref}
\usepackage{cite}
\usepackage{authblk}

\begin{document}

\title{A Brief History of Algebra with a Focus on the Distributive Law and Semiring Theory}

\author{Peyman Nasehpour}
\affil{Department of Engineering Science \\ Golpayegan University of Technology \\ Golpayegan, Isfahan Province \\ IRAN\\ \href{mailto:nasehpour@gut.ac.ir}{nasehpour@gut.ac.ir,}
	\href{mailto:nasehpour@gmail.com}{nasehpour@gmail.com}}

\maketitle
\thispagestyle{empty}

\begin{abstract}
   In this note, we investigate the history of algebra briefly. We particularly focus on the history of rings, semirings, and the distributive law.
\end{abstract}

\begin{quote}
	
	``\emph{I am sure that no subject loses more than mathematics by any attempt to dissociate it from its history.}"
	
	- J. W. L. Glaisher (1848--1928), English mathematician
	
\end{quote}

\begin{quote}
	
	``\emph{One can conceive history as an argument without end.\footnote{This is a translation of Geyl's sentence in his book with the title ``Napoleon: Voor en Tegen in de Franse Geschiedschrijving'', published in 1946. The original sentence in Dutch \cite{Geyl1946} is: ``Men kan de geschiedschrijving opvatten als een discussie zonder eind. ''}}''
	
	- Pieter Catharinus Arie Geijl (1887--1966), Dutch historian
	
\end{quote}

\section{Etymology of Algebra}

The word ``algebra" is derived from the Arabic word Al-Jabr, and this comes from the treatise written in 820 by the medieval Persian mathematician \cite{Saliba1998}, Muhammad ibn Musa al-Khwarizmi, entitled, in Arabic ``\emph{Kit\={a}b al-mukhta\d{s}ar f\={i} \d{h}isab al-\v{g}abr wa-'l-muq\={a}bala}", which can be translated as ``The Compendious Book on Calculation by Completion and Balancing" \cite{vanderWaerden1985}.\footnote{\textit{Key words and phrases}. History of Algebra, Distributive Law, Semiring Theory\\ 2010 \textit{Mathematics Subject Classification}. 01A05, 97A30}

\section{A Brief History of Algebra} 

One may consider the history of algebra to have two main stages, i.e., ``classical algebra" that mostly was devoted to solving the (polynomial) equations and ``abstract algebra", also called ``modern algebra" that is all about the study of algebraic structures, something that today algebraists do.

\subsection{Classical Algebra} 

Classical algebra can be divided into three sub-stages:

\begin{enumerate}
	\item In the very early stages of algebra, concepts of algebra were geometric, proposed by the Babylonians apparently and developed by the Greeks and later revived by Persian mathematician Omar Khayyam (1048--1131) since the main purpose of Khayyam in this matter was to solve the cubic equations with intersecting conics \cite[p. 64]{Cooke2008}.
	
	\item The main step of passing from this substage to the substage of equation-solving was taken by Persian mathematician al-Khwarizmi (c.780--c.850) in his treatise Al-Jabr that it gives a detailed account of solving polynomials up to the second degree.
	
	\item While the mathematical notion of a function was implicit in trigonometric tables, apparently the idea of a function was proposed by the Persian mathematician Sharaf al-Din al-Tusi (died 1213/4), though his approach was not very explicit, perhaps because of this point that dealing with functions without symbols is very difficult. Anyhow algebra did not decisively move to the dynamic function substage until the German mathematician Gottfried Leibniz (1646--1716) \cite{KatzBarton2007}.
	
\end{enumerate}

One may also view the history of the development of classical algebra in another perspective namely passing from rhetorical substage to full symbolic substage. According to this view, the history of algebra is divided into three substages through the development of symbolism: 1. Rhetorical algebra, 2. Syncopated algebra, and 3. Symbolic algebra \cite{Stalling2000}.\footnote{Also, see {\em The Evolution of Algebra}, Science, {\bf 18}(452) (1891), 183--187. Retrieved from \href{http://www.jstor.org/stable/1766702}{http://www.jstor.org/stable/1766702} }

\begin{enumerate}
	\item In rhetorical algebra, all equations are written in full sentences, but in symbolic algebra, full symbolism is used, the method that we do today. So far as we know, rhetorical algebra was first developed by the ancient Babylonians and was developed up to 16$^{th}$ century.
	
	\item Syncopated algebra used some symbolism but not as full as symbolic algebra. It is said that syncopated algebraic expression first appeared in the book Arithmetica by the ancient Greek mathematician Diophantus (born sometime between AD 201 and 215 and passed away sometime between AD 285 and 299) and was continued in the book Brahma Sphuta Siddhanta by the ancient Indian mathematician Brahmagupta (598--c.670 CE).
	
	\item The full symbolism can be seen in the works of the French mathematician, Ren\'{e} Descartes (1596--1650), though early steps of symbolic algebra was taken by Moroccan mathematician Ibn al-Bann\={a} al-Marak\={u}sh\={i} (1256--c.1321) and Andalusian mathematician Ab\={u} al-Hasan al-Qalas\={a}d\={i} (1412--1486) \cite[p. 162]{MerzbachBoyer2011}.
\end{enumerate}

\subsection{Modern Algebra}

The transition of algebra from the ``classical" to the ``modern" form occurred in about the middle of 19$^{th}$ century when mathematicians noticed that classical tools are not enough to solve their problems. During the time of classical algebra in the Renaissance, Italian mathematicians Scipione del Ferro (1465--1526) and Niccol\`{o} Tartaglia (1499/1500--1557) found the solution of the equations of degree 3 \cite[p. 10]{GrantKleiner2015} and another Italian mathematician Lodovico Ferrari (1522--1565) solved equations of degree 4, but, then, it was the Italian mathematician and philosopher Paolo Ruffini (1765--1822) and later the Norwegian mathematician Niels Henrik Abel (1802--1829) who used abstract algebra techniques to show that equations of degree 5 and of higher than degree 5 are not always solvable using radicals (known as Abel-Ruffini theorem \cite[\S 13]{Tignol2016}). Finally, the French mathematician \'{E}variste Galois (1811--1832) used group theory techniques to give a criterion deciding if an equation is solvable using radicals. We also need to add that the Italian mathematician Gerolamo Cardano (1501--1576), who is definitely one of the most influential mathematicians of the Renaissance, received a method of solving cubic equation from Tartaglia and promised not to publish it but he did. Since Cardano was the first to publish the explicit formula for solving the cubic equations, this is most probably why this formula is called Cardano's formula. At the end, though the British mathematician Arthur Cayley (1821--1895) was the first who gave the abstract definition of a finite group \cite{Cayley1854}, the English mathematician George Boole (1815--1864), in his book \emph{Mathematical Analysis of Logic} (1847), was most probably the first who formulated an example of a non-numerical algebra, a formal system, which can be investigated without explicit resource to their intended interpretations \cite[p. 1]{DunnHardegree2001}. It is also good to mention that the first statement of the modern definition of an abstract group was given by the German mathematician Walther Franz Anton von Dyck (1856--1934) \cite{vonDyck1882}.

\section{A Brief History of the Distributive Law} 

An algebraic structure on a set (called underlying set or carrier set) is essentially a collection of finitary operations on it \cite[p. 41, 48]{Cohn1981}. Since ring-like structures have two binary operations, often called addition and multiplication, with multiplication distributing over addition and the algebraic structure ``semiring'' is one of them, we continue this note by discussing the history of the distributive law in mathematics briefly.

The distributive law, in mathematics, is the law relating the operations of multiplication and addition, stated symbolically, $a(b + c) = ab + ac$. Ancient Greeks were aware of this law. The first six books of Elements presented the rules and techniques of plane geometry. \textrm{Book \textbf{I}} included theorems about congruent triangles, constructions using a ruler and compass, and the proof of the Pythagorean theorem about the lengths of the sides of a right triangle. \textrm{Book \textbf{II}} presented geometric versions of the distributive law $a(b + c + d) = ab + ac + ad$ and formulas about squares, such as $(a+b)^2 = a^2+2ab+b^2$ and $a^2-b^2 = (a+b)(a-b)$ \cite[p. 32]{Bradley2006Vol1}.

Apart from arithmetic, the distributive law had been noticed, years before the birth of abstract algebra, by the inventors of symbolic methods in the calculus. While the first use of the name distributive operation is generally credited to the French mathematician Fran\c{c}ois-Joseph Servois (1768--1847) \cite{Servois1814}, the Scottish mathematician Duncan Farquharson Gregory (1813--1844), who wrote a paper in 1839 entitled \emph{On the Real Nature of Symbolical Algebra}, brought out clearly the commutative and distributive laws \cite[p. 331]{Cajori1909}.

Boole in his book entitled \emph{Mathematical Analysis of Logic} (1847) mentioned this law by giving the name distributive to it and since he was, most probably, the first who formulated an example of a non-numerical algebra, one may consider him the first mathematician who used this law in abstract algebra \cite{Boole1847}. Some years later, Cayley, in a paper entitled \emph{A Memoir on the Theory of Matrices} (1858), showed that multiplication of matrices is associative and distributes over their finite addition \cite{Cayley1858}.

The distributive law appeared naturally in ring-like structures as well. Though the first axiomatic definition of a ring was given by the German mathematician Abraham Halevi (Adolf) Fraenkel (1891--1965) and he did mention the phrase ``the distributive law'' (in German ``die distributiven Gesetze"), but his axioms were stricter than those in the modern definition \cite[p. 11]{Fraenkel1914}. Actually, the German mathematician Emmy Noether (1882--1935) who proposed the first axiomatic modern definition of (commutative) rings in her paper entitled ``Idealtheorie in Ringbereichen", did also mention the phrase ``Dem distributiven Gesetz" among others such as ``Dem assoziativen Gesetz" and ``Dem kommutativen Gesetz" \cite[p. 29]{Noether1921}.

\section{A Brief History of Semirings}

The most familiar examples for semirings in classical algebra are the semiring of nonnegative integers or the semiring of nonnegative real numbers. The first examples of semirings in modern algebra appeared in the works of the German mathematician Richard Dedekind (1831--1916) \cite{Dedekind1894}, when he worked on the algebra of the ideals of rings \cite{Golan1999(b)}. In point of fact, it was Dedekind who proposed the concept of ideals in his earlier works on number theory, as a generalization of the concept of ``ideal numbers" developed by the German mathematician Ernst Kummer (1810--1893) \cite{Kummer1847}. Others such as the English mathematician Francis Sowerby Macaulay (1862--1937) and the German mathematicians Emanuel Lasker (1868--1941), Emmy Noether (1882--1935), Wolfgang Krull (1899--1971) and Paul Lorenzen (1915--1994) also studied the algebra of ideals of rings \cite{Krull1924, Lasker1905, Lorenzen1939, Macaulay1916, Noether1921}. Semirings appeared implicitly in the works of the German mathematician David Hilbert (1862--1943) and the American mathematician Edward Vermilye Huntington (1874--1952) in connection with the axiomatization of the natural and nonnegative rational numbers \cite{Golan1999(b), Hilbert1899, Huntington1902}. But, then, it was the American mathematician Harry Schultz Vandiver (1882--1973) who used the term ``semi-ring" in his 1934 paper entitled ``Note on a simple type of algebra in which cancellation law of addition does not hold" for introducing an algebraic structure with two operations of addition and multiplication such that multiplication distributes on addition, while cancellation law of addition does not hold \cite{Vandiver1934}. The foundations of algebraic theory for semirings were laid by Samuel Bourne and others in the l950's \cite[p. 6]{AhsanMordesonShabir2012}. For example, the concept of ideals for semirings was introduced by Samuel Bourne \cite{Bourne1951}. In the years between 1939 and 1956, Vandiver published at least six more papers on semirings \cite{Vandiver1939,Vandiver1940,Vandiver1952,Vandiver1954,VandiverWeaver1955,VandiverWeaver1956}, but it seems he was not successful to draw the attention of mathematicians to consider semirings as an independent algebraic structure that is worth to be developed \cite{Golan2005}. In fact, semirings, as most of the other concepts in mathematics, were not developed as an exercise for generalization, only for the sake of generalization! Actually, in the late 1960s, semirings were considered as a more serious topic by researchers when real applications were found for them. The Polish-born American mathematician Samuel Eilenberg (1913--1998) and a couple of other mathematicians started developing formal languages and automata theory systematically \cite{Eilenberg1974}, which they have strong connections to semirings. Since then many mathematicians and computer scientists have broadened the theory of semirings and related structures \cite{Glazek2002}. Definitely, the reference books on semirings and other ring-like algebraic structures including the books \cite{AhsanMordesonShabir2012,Bistarelli2004, DrosteKuichVogler2009, Golan1999(a), Golan1999(b), Golan2003, HebischWeinert1998, KuichSalomaa1986, Pilz1983} have helped to the popularity of these rather new algebraic structures. Today many journals specializing in algebra have editors who are responsible for semirings. Semirings not only have significant applications in different fields such as automata theory in theoretical computer science, (combinatorial) optimization theory, and generalized fuzzy computation, but are fairly interesting generalizations of two broadly studied algebraic structures, i.e., rings and bounded distributive lattices \cite{Glazek2002,Golan1999(b)}. The number of publications in the field of semiring theory, the beauty of the semirings and their broad applications in different areas of science should convince us that today semiring theory is an established one and its development, even in pure mathematics, is valuable and important.

\section{References for further studies} In order to have a bit more of ``an argument without end", the reader may like to refer the books \cite{Bradley2006Vol1, Bradley2006Vol2, Bradley2006Vol3, Bradley2006Vol4, Bradley2006Vol5,Wardhaugh2010} on the general history of mathematics and the books \cite{Alten2014,Kleiner2007,vanderWaerden1985} on the history of algebra.

\section*{Acknowledgements} The author's main interests in algebra are in commutative algebra and semiring theory. It is a pleasure to thank both Professor Dara Moazzami and Professor Winfried Bruns to help and encourage the author in order to work on these fields of algebra. The author is supported in part by the Department of Engineering Science at Golpayegan University of Technology and his special thanks go to the Department for providing all the necessary facilities available to him for successfully conducting this research.

\bibliographystyle{plain}

\end{document}